\documentclass[11pt,a4paper]{amsart}
\usepackage[english]{babel}
\usepackage{graphicx}
\usepackage{framed}
\usepackage[normalem]{ulem}
\usepackage{amsmath}
\usepackage{amsthm}
\usepackage{amssymb}
\usepackage{amsfonts}
\usepackage[foot]{amsaddr}
\usepackage{mathrsfs,amsmath}
\usepackage{enumerate}
\usepackage[utf8]{inputenc}
\usepackage[top=1 in,bottom=1in, left=1 in, right=1 in]{geometry}
\usepackage{mathtools}
\usepackage{tikz-cd}
\usepackage[font=small,labelfont=bf]{caption} 
\usepackage[subrefformat=parens]{subcaption}
\usepackage{cite}
\usepackage{abstract}

\usetikzlibrary{decorations.markings, arrows.meta}

\tikzset{
  marrow/.style={decoration={markings,mark=at position 0.5 with {\arrow{#1}}}, postaction=decorate}
}

\usepackage{enumitem}
\setlist{  
  listparindent=\parindent,
  parsep=0pt,
}

\usepackage{bbm}

\usepackage{stackrel}

\newcommand{\R}{\mathbb{R}}
\newcommand{\N}{\mathbb{N}}

\newtheorem{thm}{Theorem}

\newtheorem{prop}[thm]{Proposition}

\theoremstyle{definition}

\newtheorem{conj}[thm]{Conjecture}

\theoremstyle{remark}

\newtheorem{rmk}[thm]{Remark}

\newcommand{\eproof}{\hfill\qed}

\newcommand{\bproofof}[1]{\noindent{\textit{Proof of #1. }}}

\begin{document}

\title{Legendrians with vanishing Shelukhin-Chekanov-Hofer metric}
\author{Lukas Nakamura}
\maketitle

\begin{abstract}
	We show that the Legendrian lift of an exact, displaceable Lagrangian has vanishing Shelukhin-Chekanov-Hofer pseudo-metric by lifting an argument due to Sikorav to the contactization. In particular, this proves the existence of such Legendrians, providing counterexamples to a conjecture of Rosen and Zhang.\\
\end{abstract}

Let $(M^{2n+1}, \alpha)$ be a strict contact manifold. For two Legendrians $L_0$ and $L_1$ which are Legendrian isotopic, denote by $d_\alpha (L_0,L_1)$ the Shelukhin-Chekanov-Hofer (pseudo-)metric induced by $\Vert H_t \Vert \coloneqq \int_0^1 \underset{x \in M}{\max} |H_t(x)| dt$ via $d_\alpha (L_0,L_1) \coloneqq \underset{H_t}{\inf} \Vert H_t \Vert$, where the infimum is taken over all compactly supported Hamiltonians $H_t: M \to \R$ whose associated contact isotopy $\phi_t$ satisfies $\phi_1(L_0) = L_1$.

Recall the following result by Rosen and Zhang.

\begin{thm}\label{thm:chekanovs dichotomy for legendrians}\hspace{-1mm}\textnormal{\cite{rz18}}\,
	Let $L \subseteq M$ be a properly embedded Legendrian submanifold. Then on the Legendrian isotopy class of $L$, $d_\alpha$ is either non-degenerate or vanishes identically. 
\end{thm}

They conjectured that for closed Legendrians, $d_\alpha$ is always non-degenerate. We show that this is not the case by lifting an example due to Sikorav of an exact Lagrangian with vanishing Chekanov-Hofer metric $d$ to the contactization.

The following proposition is a reformulation of (a generalization of) Sikorav's example in \cite{che00}.

\begin{prop}\label{prop:lagrangian with vanishing metric}
	Let $(M,d \lambda)$ be an exact symplectic manifold with complete Liouville flow, and let $L_0,L_1 \subseteq M$ be two closed, exact Lagrangian submanifolds which are Hamiltonian isotopic and disjoint. Then there exists a 1-parameter\footnote{$H^s_t$ depends continuously on $s$ with respect to the $C^\infty$-topology.} family $H^s_t:M \to \R, t \in [0,1], s \in (0,\infty),$ of compactly supported Hamiltonians with $\Vert H^s_t \Vert = s$ so that the associated Hamiltonian isotopies $\phi^{s}_t$ satisfy $\phi^{s}_1(L_0)=L_1$ and $\phi^{s}_1|_{L_0} = \phi^{1}_1|_{L_0}$ for all $s$.
\end{prop}

\begin{rmk}
	In particular, this shows that $d(L_0,L_1)=0$, which implies that the Chekanov-Hofer pseudo-metric vanishes for $L_0$ by Chekanov's dichotomy \cite[Theorem 2]{che00}. 
\end{rmk}

\begin{rmk}\label{rmk:exact displaceable Lagrangians}
	One example of Lagrangians satisfying the conditions in the theorem is constructed in \cite{mul90}. More recently, Murphy showed in \cite{mur13} that the isotopy class of any closed, formal Lagrangian embedding into the symplectization of an overtwisted contact manifold of dimension at least $5$ can be realized by an exact Lagrangian. Such a Lagrangian can always be displaced by shifting it in the symplectization direction via the Liouville flow, and such Lagrangian isotopies are induced by Hamiltonian isotopies by exactness of the Lagrangian.
\end{rmk}

\bproofof{Proposition \ref{prop:lagrangian with vanishing metric}}
Let $Z$ denote the Liouville vector field on $M$ defined by $i_Z d \lambda = \lambda$, and let $\Phi_t$ denote the corresponding Liouville flow. Recall that the Liouville flow satisfies $\Phi_t^*\lambda = e^t \lambda$. This implies that $\Phi_t(L_i)$ is an exact Lagrangian for all $t \in \R$ and $i \in \{0,1\}$. Therefore, the homotopy $\Phi_t(L_0 \coprod L_1)$ of Lagrangian embeddings is exact, and thus there exists a Hamiltonian isotopy $\Psi_t$ on $M$ with $\Psi_t|_{L_0 \coprod L_1} = \Phi_t|_{L_0 \coprod L_1}$ for all $t \in \R$. Let $H_t:M \to \R$ be a compactly supported, time-dependent Hamiltonian with associated Hamiltonian isotopy $\phi_t$ such that $\phi_1(L_0)=L_1$. Let $s \in \R$ be an arbitrary number. Then the Hamiltonian $H_t \circ \Psi_s^{-1}$ generates the isotopy $\Psi_s \phi_t \Psi_s^{-1}$, and the Hamiltonian $H^s_t \coloneqq e^{-s} H_t \circ \Psi_s^{-1} \circ \Phi_s$ generates the Hamiltonian isotopy $\Phi_s^{-1} \Psi_s \phi_t \Psi_s^{-1} \Phi_s$ which satisfies $\Phi_s^{-1} \Psi_s \phi_1 \Psi_s^{-1} \Phi_s|_{L_0} = \phi_1|_{L_0}$  by definition of $\Psi_s$. Note that $\Vert H^s_t \Vert = e^{-s} \Vert H_t \Vert$. After a reparametrization of the family $\{H^s_t\}_s$, we find the desired family of Hamiltonians.
\eproof\\

By lifting the Hamiltonians in Proposition \ref{prop:lagrangian with vanishing metric} to the contactization we prove the following.

\begin{thm}\label{thm:legendrian lift of closed, exact lagrangian has vanishing d_alpha}
	Under the assumptions of Proposition \ref{prop:lagrangian with vanishing metric}, the Legendrian lift $\Lambda_0$ of $L_0$ to the contactization $(M \times \R,dz + \lambda)$ has vanishing Shelukhin-Chekanov-Hofer metric.
\end{thm}

We expect the following generalization to hold as well.

\begin{conj}
	\textit{Let $(M,d \lambda)$ be an exact symplectic manifold and $L$ be a closed exact Lagrangian submanifold with vanishing Chekanov-Hofer metric. Then its Legendrian lift $\Lambda$ has vanishing Shelukhin-Chekanov-Hofer metric.}
\end{conj}

\begin{rmk}
	After completion of this manuscript, it came to our attention that Cant \cite[Theorem 2]{can23} independently also proved the existence of Legendrians with vanishing Shelukhin-Chekanov-Hofer distance. To be more precise, he showed that for two contact manifolds $Y_0$ and $Y_1$ and an exact Lagrangian submanifold $L$ in the symplectization $SY_1$ of $Y_1$, there exist Legendrian submanifolds of $Y_0 \times S Y_1$ with vanishing $d_\alpha$. The construction in \cite{can23} is similar to our construction, but the proof of the degeneracy of $d_\alpha$ in \cite{can23} differs from our proof presented below.
\end{rmk}

\bproofof{Theorem \ref{thm:legendrian lift of closed, exact lagrangian has vanishing d_alpha}}
Let $H^s_t:M \to \R$ be a Hamiltonian as in Proposition \ref{prop:lagrangian with vanishing metric}. Let $G^s_t(x,z) \coloneqq H^s_t(x)$ be the cylindrical lift of $H^s_t$ to $M \times \R$ with associated contact vector field $Y^s_t = (X^s_t,f^s_t \partial_z)$. The defining equations for $Y^s_t$ are $H^s_t = f^s_t + \lambda(X^s_t)$ and $\pi^*(\iota_{X^s_t} d \lambda)|_{\ker(dz+\lambda)} = -d G^s_t|_{\ker(dz + \lambda)}$, where $\pi:M\times \R \to M$ denotes the projection onto $M$. The latter equation is equivalent to $\iota_{X^s_t} d \lambda = -d H^s_t$ which means that $X^s_t$ is the Hamiltonian vector field associated to $H^s_t$. This implies that the contact isotopy $\psi^s_t$ associated to $G^s_t$ is of the form $\psi^s_t(x,z) = (\phi^s_t(x), \rho^s_t(x) + z)$, where $\phi^s_t$ denotes the Hamiltonian isotopy associated to $H^s_t$. Let $\Lambda_0$ be a Legendrian lift of $L_0$ to $M \times \R$, ie a Legendrian submanifold of $(M \times \R,dz + \lambda)$ which is a lift of $L_0$ with respect to the canonical projection $M \times \R \to M$. Then $\psi^s_1(\Lambda_0)$ is a Legendrian lift of $L_1$. Since the Legendrian lift is unique up to a shift in the $z$-direction and $\phi^s_1|_{L_0} = \phi^1_1|_{L_0}$, it follows that $g^s \coloneqq \rho^s_1(x) - \rho^1_1(x), x \in L_0,$ does not depend on the choice of $x \in L_0$. Furthermore, $g^s$ depends continuously on $s$. After cutting off $G^s_t$ outside of a sufficiently large (possibly $s$-dependent) compact set, we may assume that $\psi^s_t$ has compact support, $\Vert G^s_t \Vert = s$, and $\psi^s_1(x,z)=(\phi^1_1(x),g^s + \rho^1_1(x) + z)$ for all $(x,z) \in \Lambda_0$, $s \in (0,\infty)$. In particular, $\phi^R_{g^{s_1}-g^{s_0}} \circ \psi^{s_0}_1|_{\Lambda_0} = \psi^{s_1}_1|_{\Lambda_0}$ for all $s_0,s_1 \in (0,\infty)$, where $\phi^R_t(x,z)=(x,z+t)$ denotes the time-$t$ Reeb flow.

First, we assume that $g^s$ has a convergent subsequence as $s \to 0$, i.e. there exists a sequence $\{s_i\}_{i \in \N}$ with $s_i \to 0$ and $g^{s_i} \to g^0 \in \R$ as $i \to \infty$. Then the image $\Lambda_\infty$ of $x \mapsto (\phi^1_1(x),g^0 + \rho^1_1(x)), x \in L_0,$ is a Legendrian lift of $L_1$, and it satisfies 
\begin{equation}
d_{dz+\lambda}(\Lambda_\infty,\Lambda_0) \leq d_{dz+\lambda}(\Lambda_\infty,\psi^{s_i}_1(\Lambda_0)) + d_{dz+\lambda}(\psi^{s_i}_1(\Lambda_0), \Lambda_0) \leq |g^0-g^{s_i}| + s_i \to 0, \quad i \to \infty,
\end{equation}
where $d_{dz+\lambda}(\Lambda_\infty,\psi^{s_i}_1(\Lambda_0)) \leq |g^0-g^{s_i}|$ because $\Lambda_\infty$ and $\psi^{s_i}_1(\Lambda_0)$ are related by a $(g^0-g^{s_i})$-shift in the Reeb direction.

In the case that $g^s$ has no convergent subsequence as $s \to 0$, it follows that $g^s \to \pm \infty$ as $s \to 0$ by continuity of $s \mapsto g^s$. Again by continuity of $s \mapsto g^s$, we can then find sequences $\{s_i\}_{i\in\N}$ and $\{s'_i\}_{i \in \N}$ with $s_i, s'_i \to 0$ as $i \to \infty$ and $g^{s_i}-g^{s'_i}=1$ for all $i \in \N$. It follows that 
\begin{equation}
\begin{split}
d_{dz+\lambda}(\psi^{s_1}(\Lambda_0),\psi^{s'_1}(\Lambda_0)) &= d_{dz+\lambda}(\psi^{s_i}(\Lambda_0),\psi^{s'_i}(\Lambda_0))\\
&\leq d_{dz+\lambda}(\psi^{s_i}(\Lambda_0),\Lambda_0) + d_{dz+\lambda}(\Lambda_0,\psi^{s'_i}(\Lambda_0)) \to 0, \quad i \to \infty,
\end{split}
\end{equation}
where the equality in the first line is due to the fact that the pairs $(\psi^{s_1}(\Lambda_0),\psi^{s'_1}(\Lambda_0))$ and $(\psi^{s_i}(\Lambda_0),\psi^{s'_i}(\Lambda_0))$ are related by the time-$(g^{s_i}-g^{s_1})$ Reeb flow, which is a strict contact isotopy.

In either case, we see that the Shelukhin-Chekanov-Hofer metric on the Legendrian isotopy class of $\Lambda_0$ is degenerate and hence vanishes identically by Theorem \ref{thm:chekanovs dichotomy for legendrians}.
\eproof \\

\textbf{Open Questions.} In Murphy's construction of displaceable exact Lagrangians, the contactization of the ambient symplectic manifold is always overtwisted. This observation leads to the following question:\\

\textbf{Question 1:} \textit{Are there tight contact manifolds which admit closed Legendrian submanifolds with vanishing Shelukhin-Chekanov-Hofer metric?}\\

Related to this question, it is not clear whether non-degeneracy of the Shelukhin-Chekanov-Hofer metric for Legendrians actually is a property of Legendrian isotopy classes of Legendrian embeddings, or rather a property of the ambient contact manifold. In other words:\\

\textbf{Question 2:} \textit{Are there contact manifolds which simultaneously admit 1) an isotopy class of closed Legendrian submanifolds with non-degenerate $d_\alpha$ and 2) an isotopy class of closed Legendrian submanifolds with vanishing $d_\alpha$?}\\

\textbf{Acknowledgments:} I am grateful to my adviser Tobias Ekholm for his guidance and interesting discussions about the contents of this paper. I would also like to thank Georgios Dimitroglou Rizell for directing my attention to Murphy's work on the existence of closed and exact Lagrangians. I am also grateful to the anonymous referee whose comments helped to improve this article.

This work is partially supported by the Knut and Alice Wallenberg Foundation, grant KAW2020.0307.\\

\bibliography{references}
\bibliographystyle{amsalpha}

\end{document}